\documentclass{article}

\usepackage{
 amsmath,
 amsxtra,
 amsthm,
 amssymb,
 etex,
 mathrsfs,
 mathtools,
 }
\usepackage[all]{xy}
\usepackage[T2A,T1]{fontenc}
\usepackage[OT2,T1]{fontenc}
\newcommand{\Zhe}{\mbox{\usefont{T2A}{\rmdefault}{m}{n}\CYRZH}}

\newtheorem{theorem}{Theorem}[section]
\newtheorem{lemma}[theorem]{Lemma}
\newtheorem{conjecture}[theorem]{Conjecture}
\newtheorem{proposition}[theorem]{Proposition}
\newtheorem{corollary}[theorem]{Corollary}

\theoremstyle{definition}
\newtheorem{defn}[theorem]{Definition}
\newtheorem{remark}[theorem]{Remark}

\newtheorem*{acknowledgments}{Acknowledgments}

\newcommand{\bd}{\begin{defn}}
\newcommand{\ed}{\end{defn}}
\newcommand{\bl}{\begin{lemma}}
\newcommand{\el}{\end{lemma}}
\newcommand{\bp}{\begin{proposition}}
\newcommand{\ep}{\end{proposition}}
\newcommand{\bt}{\begin{theorem}}
\newcommand{\et}{\end{theorem}}
\newcommand{\bc}{\begin{corollary}}
\newcommand{\ec}{\end{corollary}}
\newcommand{\br}{\begin{remark}}
\newcommand{\er}{\end{remark}}
\newcommand{\ba}{\begin{array}}
\newcommand{\ea}{\end{array}}
\newcommand{\bpf}{\begin{proof}}
\newcommand{\epf}{\end{proof}}

\newcommand{\Z}{\mathbb{Z}}
\newcommand{\Q}{\mathbb{Q}}
\newcommand{\Zp}{\mathbb{Z}_{p}}
\newcommand{\Qp}{\mathbb{Q}_{p}}

\newcommand{\Ep}{E_{p^{\infty}}}
\newcommand{\al}{\alpha}
\newcommand{\be}{\beta}

\newcommand{\ga}{\gamma}

\newcommand{\La}{\Lambda}

\DeclareMathOperator{\Sel}{Sel} \DeclareMathOperator{\Gal}{Gal}
 \DeclareMathOperator{\rank}{rank}

\DeclareMathOperator{\Ext}{Ext}

\newcommand{\cyc}{\mathrm{cyc}}

\newcommand{\M}{\mathfrak{M}}

\newcommand{\mM}{\mathcal{M}}

\newcommand{\ot}{\otimes}
\newcommand{\ilim}{\displaystyle \mathop{\varinjlim}\limits}

\newcommand{\coker}{\mathrm{coker}\,}

\newcommand{\lra}{\longrightarrow}

\newcommand{\ps}[1]{[[ #1 ]]}

  \DeclareFontFamily{U}{wncy}{}
  \DeclareFontShape{U}{wncy}{m}{n}{<->wncyr10}{}
  \DeclareSymbolFont{mcy}{U}{wncy}{m}{n}
  \DeclareMathSymbol{\sha}{\mathord}{mcy}{"58}

\begin{document}
\title{On pseudo-nullity of fine Mordell-Weil group}
 \author{
  Meng Fai Lim\footnote{School of Mathematics and Statistics,
Central China Normal University, Wuhan, 430079, P.R.China.
 E-mail: \texttt{limmf@ccnu.edu.cn}} \quad
 Chao Qin\footnote{College of Mathematical Sciences,
Harbin Engineering University,
Harbin, 150001, P.R.China.
 E-mail: \texttt{qinchao@hrbeu.edu.cn}} 
\quad
 Jun Wang\footnote{Institute for Advanced Study in Mathematics of HIT,
Harbin Insitute of Technology,
Harbin, 150001, P.R.China.
 E-mail: \texttt{junwangmath@hit.edu.cn}} }
\date{}
\maketitle

\begin{abstract} \footnotesize
\noindent
Let $E$ be an elliptic curve defined over $\Q$ with good ordinary reduction at a prime $p\geq 5$, and let $F$ be an imaginary quadratic field. Under appropriate assumptions, we show that the Pontryagin dual of the fine Mordell-Weil group of $E$ over the $\Zp^2$-extension of $F$ is pseudo-null as a module over the Iwasawa algebra of the group $\Zp^2$.


\end{abstract}

\section{Introduction}
At the turn of the millennium, Coates and Sujatha \cite{CS05}, and a little later, Wuthrich \cite{Wu} initiated a systematic study on the fine Selmer group of an elliptic curve $E$. The fine Selmer group is a subgroup of the classical $p$-primary Selmer group defined by stricter local conditions at primes above $p$. Analogous to the usual Selmer group, this fine Selmer group $R(E/F)$ sits in the
middle of the following short exact sequence
\[0\lra \mM(E/F) \lra R(E/F) \lra \Zhe(E/F) \lra 0,\]
where $\mM(E/F)$ and $\Zhe(E/F)$ are the fine Mordell-Weil group and fine Tate-Shafarevich group respectively (defined in the sense of Wuthrich \cite{WuTS}) which can be thought as the ``fine'' counterpart of the usual Mordell-Weil group and Tate-Shafarevich group.

In \cite[Conjecture B]{CS05}, Coates-Sujatha proposed the following conjecture.

\begin{conjecture}[Conjecture B] Let $E$ be an elliptic curve defined over a number field $F$. Suppose that $F_\infty$ is a $p$-adic Lie extension of $F$ for which $\Gal(F_\infty/F)$ has dimension $\geq 2$ containing the cyclotomic $\Zp$-extension $F^\cyc$. Then $R(A/F_\infty)^\vee$ is pseudo-null over $\Zp\ps{\Gal(F_\infty/F)}$.
\end{conjecture}

Here a finitely generated $\Zp\ps{\Gal(F_\infty/F)}$-module $M$ is said to be pseudo-null if \[\Ext^i_{\Zp\ps{\Gal(F_\infty/F)}}(M, \Zp\ps{\Gal(F_\infty/F)}) =0\]
 for $i=0,1$.

This conjecture remains wide open but we refer readers to \cite{Bh, Jh, LimPS,LimFine} for some discussion and numerical examples in support of the conjecture. 
The goal of this note is to provide further theoretical support to Conjecture B. Namely, we will prove the following.

\bt \label{mainthm}
 Let $E$ be an elliptic curve defined over $\Q$ with good ordinary reduction at a prime $p\geq 5$, and let $F$ be an imaginary quadratic field. Suppose that $E$ has no complex multiplication, and that the discriminant of $F$ is coprime to the conductor of $E$. In the event that the root number $\epsilon(E/F, 1)$ equals $-1$, assume further that $p$ does not divide the class number of $F$. Suppose that $\Sel(E/F^\cyc)$ is cofinitely generated over $\Zp$. Then the Pontryagin dual $\mathcal{M}(E/F_\infty)^\vee$ of the fine Mordell-Weil group is a pseudo-null $\Zp\ps{\Gal(F_\infty/F)}$-module, where $F_\infty$ is the $\Zp^2$-extension of $F$.  
\et

As an application, we can prove the following, where we note that every $\Zp$-extension of $F$ is necessarily contained in $F_\infty$.

\bc \label{maincor1}
Retain the settings of Theorem \ref{mainthm}. Then for every $\Zp$-extension $L_\infty$ of $F$, we have that $\mathcal{M}(E/L_\infty)^\vee$ is torsion over $\Zp\ps{\Gal(L_\infty/L)}$.
\ec

It is conjectured that $R(E/L_\infty)^\vee$ should always be torsion over $\Zp\ps{\Gal(L_\infty/L)}$ (see \cite{PR00, WuPhD, LimFineDoc}). The above corollary thus constitutes a partial evidence towards this conjecture. Furthermore, if one is willing to assume the finiteness of the fine Tate-Shafarevich group, we have the following observation.

\bc \label{maincor2}
Retain the settings of Theorem \ref{mainthm}. Suppose that $L_\infty$ is a $\Zp$-extension of $F$  with the property that $\Zhe(E/L)$ is finite for every finite extension $L$ of $F$ contained in $L_\infty$. Then $R(E/L_\infty)^\vee$ is torsion over $\Zp\ps{\Gal(L_\infty/L)}$.
\ec

The proofs of the preceding theorem and corollaries will be given in Section \ref{fine selmer sec}.

\section{Selmer group}

We now introduce the Selmer group in a slightly more general context. As a start, we let $F$ be an arbitrarily fixed number field and $E$ an elliptic curve defined over $F$. Let $S$ be a finite set of primes of $F$ which contains all the primes above $p$, the infinite primes and the primes of bad reduction of $E$. We shall also write $S_p$ for the set of primes of $F$ above $p$. Denote by $F_S$ the maximal algebraic extension of $F$ which is unramified outside $S$. If $\mathcal{L}$ is a (possibly infinite) extension of $F$ contained in $F_S$ and $S'\subseteq S$, we write $S'(\mathcal{L})$ for the set of primes of $\mathcal{L}$ above $S'$.

For each $v\in S$ and a finite extension $L$ of $F$, set
\[ J_v(E/L) = \bigoplus_{w|v}H^1(L_w,E)_{p^\infty}.\]
If $\mathcal{L}$ is an infinite extension of $F$ contained in $F_S$, we define
\[ J_v(E/\mathcal{L}) = \ilim_L J_v(E/L),\]
where $L$ runs through all finite extensions  of $F$ contained in $\mathcal{L}$.

The classical ($p$-primary) Selmer group of $E$ over $\mathcal{L}$ is defined by
\[ \Sel(E/\mathcal{L})=\ker\Big(H^1(G_S(\mathcal{L}), \Ep)\lra \bigoplus_{v\in S}J_v(E/\mathcal{L})\Big),\]
where we write $G_S(\mathcal{L}) = \Gal(F_S/\mathcal{L})$. The Pontryagin dual of $\Sel(E/\mathcal{L})$ is then denoted by $X(E/\mathcal{L})$.

From now on, we will always assume that the elliptic curve $E$ has good ordinary reduction at all primes of $F$ above $p$.
For each prime $v$ of $F$ above $p$, denote by $\hat{E}_v$ and $\tilde{E}_v$ the formal group of $E$ at $v$ and the reduced curve of $E$ at $v$ respectively. Furthermore, we have a short exact sequence
\[ 0\lra \hat{E}_{v, p^\infty}\lra \Ep\lra \tilde{E}_{v,p^\infty}\lra 0 \]
of discrete $\Gal(\bar{F}_v/F_v)$-modules. Since $E$ is assumed to have good ordinary reduction, both $\hat{E}_{v, p^\infty}$ and $\tilde{E}_{v,p^\infty}$ are cofree $\Zp$-modules of corank 1.

For our purposes, it is convenient to work with an equivalent description of the local terms $J_v(E/\mathcal{L})$, following an insight of Coates-Greenberg \cite{CG}.
Let $\mathcal{L}$ be an algebraic extension of $F$. For every non-archimedean prime $w$ of $\mathcal{L}$, write $\mathcal{L}_w$ for the union of the completions at $w$ of the finite extensions of $F$ contained in $\mathcal{L}$. If $w$ is a prime above $p$, we write $E_w = E_v$, where $v$ is a prime of $F$ below $w$. Finally, we shall always denote by $F^{\cyc}$ the cyclotomic $\Zp$-extension of $F$. With these in hand, we have the following lemma.

\bl \label{local coh description}
Let $\mathcal{L}$ be an algebraic extension of $F^\cyc$ which is unramified outside a set of finite primes of $F$. Then we have an isomorphism
\[ J_v(E/\mathcal{L})\cong \begin{cases} \ilim_{\mathcal{L}'}\bigoplus_{w|v}H^1(\mathcal{L}'_w, \tilde{E}_{w,p^\infty}),  & \mbox{if $v$ divides $p$}, \\
\ilim_{\mathcal{L}'}\bigoplus_{w|v}H^1(\mathcal{L}'_w, \Ep), & \mbox{if $v$ does not divide $p$},\end{cases} \]
where the direct limit is taken over all finite extensions $\mathcal{L}'$ of $F^\cyc$ contained in $\mathcal{L}$.
\el

\bpf
See \cite[Propositions 4.1, 4.7 and 4.8]{CG} or \cite[Lemma 4.1]{LimMHG}.
\epf

We now establish a control theorem for the Selmer group over a $\Zp^2$-extension.

\bp \label{cyclotomic control}
Let $E$ be an elliptic curve defined over $F$ which has good ordinary reduction at each prime of $F$ above $p$. Suppose that $F_\infty$ is a $\Zp^2$-extension of $F$ which contains $F^\cyc$. Write $H=\Gal(F_\infty/F^\cyc)$. Then the restriction map
\[   \Sel(E/F^\cyc)\lra \Sel(E/F_\infty)^H\]
has finite kernel and cokernel.
\ep

\bpf
The proposition should definitely be well-known. For the convenience of the readers, we shall supply an argument here. 
 Consider the following commutative diagram
 \[   \entrymodifiers={!! <0pt, .8ex>+} \SelectTips{eu}{}\xymatrix{
    0 \ar[r] &\Sel(E/F^\cyc)  \ar[d]^{} \ar[r] &  H^1\big(G_{S}(F^\cyc),\Ep\big)
    \ar[d]^{h} \ar[r] & \displaystyle \bigoplus_{v\in S}J_v(E/F^\cyc) \ar[d]^{\oplus g_v} \ar[r]& 0 \\
    0 \ar[r]^{} & \Sel(E/F_\infty)^{H} \ar[r]^{} & H^1\big(G_{S}(F_\infty),\Ep\big)^{H} \ar[r] &
    \left(\displaystyle \bigoplus_{v\in S}J_v(E/F_\infty)\right)^{H} &  } \]
    with exact rows. Via the snake lemma, it suffices to show that $\ker h$, $\ker g_v$ and $\coker h$ are finite. To begin with, we show that $h$ is surjective with a finite kernel. Indeed, since $H\cong \Zp$, the
restriction-inflation sequence tells us that the map $h$ is surjective with kernel $H^1(H, E(F_\infty)_{p^\infty})$. On the other hand, we have
\[ 0 = \rank_{\Zp\ps{H}}\Big( E(F_{\infty})_{p^\infty}\Big)^\vee = \rank_{\Zp} \Big(\big(E(F_{\infty})_{p^\infty}\big)^\vee\Big)_H- \rank_{\Zp} \Big(\big(E(F_{\infty})_{p^\infty}\big)^\vee\Big)^H\]
where the second equality follows from \cite[Proposition 5.3.20]{NSW}. But observe that
\[\Big(\big(E(F_{\infty})_{p^\infty}\big)^\vee\Big)_H = \Big(\big(E(F_{\infty})_{p^\infty}\big)^H\Big)^\vee = \Big(E(F^\cyc)_{p^\infty}\Big)^\vee,  \]
and the latter is finite by a theorem of Imai \cite{Imai}. Hence $\Big(\big(E(F_{\infty})_{p^\infty}\big)^\vee\Big)^H$ is also finite. But this group is precisely
$H^1(H, E(F_\infty)_{p^\infty})^\vee$, and so we have our claim.

It remains to show that $g_v$ has finite kernel for every $v$. Again, by the restriction-inflation sequence, we have $\ker g_v = \oplus_{w|v} H^1(H_w, D)$, where here the sum runs over all the primes of $F^\cyc$ above $v$, $H_w$ is the decomposition group of $w$ in $H$ and $D$ denotes either $\tilde{E}_{v,p^\infty}$ or $\Ep$ according as $v$ divides $p$ or not.
In particular, if $H_w=1$, then $H^1(H_w, D)=0$. This is indeed the case when $v$ does not divide $p$ for a $\Zp^2$-extension is unramified outside $p$. It remains to consider primes $w$ which divides $p$ and that $H_w$ is nontrivial. Since $H\cong\Zp$, it then follows that $H_w\cong \Zp$. We may then apply the same argument as in the preceding paragraph to conclude that $H^1(H_w, \tilde{E}_{v,p^\infty})$ is finite. The proof of the proposition is now complete.
\epf

The following is a corollary of the preceding proposition.

\bc \label{cyclotomic control coro}
Retain the setting as in Proposition \ref{cyclotomic control}. Then $\Sel(E/F^\cyc)$ is cofinitely generated over $\Zp$ if and only if $\Sel(E/F_\infty)$ is cofinitely generated over $\Zp\ps{H}$.
\ec

\section{Mordell-Weil group over $\Zp^2$-extension}

We continue to suppose that $F_\infty$ is a $\Zp^2$-extension of $F$ which contains $F^\cyc$. Write $H=\Gal(F_\infty/F^\cyc)$ and denote by $H_n$ the unique subgroup of $H$ with index $p^n$. The fixed field of $H_n$ is in turn denoted by $K_n$. The following hypothesis will be in full force for the remainder of the section.

\medskip
$\mathbf{(Fg)}$ $\Sel(E/F^\cyc)$ is cofinitely generated over $\Zp$.

\medskip
Write $\Lambda=\Zp\ps{H}$. We shall identify this latter ring with the power series ring $\Zp\ps{T}$ in one variable. By abuse of notation, we shall also write $\Lambda$ for the ring $\Zp\ps{T}$. Denote by $\Phi_n$ the $p^n$th-cyclotomic polynomial which is viewed as an element in $\Lambda$.

\bp \label{MW control}
Suppose that $\mathbf{(Fg)}$ is valid. Write $\La=\Zp\ps{H}$. Then the following statements are valid.

\begin{enumerate}
\item[$(a)$] There is an injective $\La$-homomorphism
\[ \big(E(F_\infty)\ot\Qp/\Zp\big)^\vee \lra \La^{\oplus r}\oplus \Big(\bigoplus_{n\geq 0}(\La/\Phi_n)^{\oplus t_n}\Big)\]
with a finite cokernel, where $\{t_n\}$ is a sequence of nonnegative integers with $t_n = 0$ for $n\gg0$.

\item[$(b)$] Define
\[ e_n = \left\{
           \begin{array}{ll}
            \displaystyle\frac{\rank_{\Z}E(K_n) - \rank_{\Z}E(K_{n-1}) }{p^n-p^{n-1}} , & \mbox{if $n\geq 1;$} \\
            \\
            \rank_{\Z}E(F^\cyc) , & \mbox{if $n=0$.}
           \end{array}
         \right.
\]
Then we have
\[ T_p\big(E(K_n)\ot \Qp/\Zp \big) \cong \bigoplus_{j=0}^n(\La/\Phi_j)^{\oplus e_j}.\]

\end{enumerate}
\ep

\bpf
By hypothesis $\mathbf{(Fg)}$, we have that $\Sel(E/K_n)$ is cofinitely generated over $\Zp$ for every $n$ (for instance, see \cite[Corollary 3.4]{HM}). Thus, it follows that $E(K_n)\ot\Qp/\Zp$ is also cofinitely generated over $\Zp$ for every $n$. In view of this observation, we may apply a similar argument to that of Lee in \cite[Theorem 2.1.2]{Lee} to obtain the conclusion of statement (a). Once, we have (a), statement (b) will follow by carrying out the argument in \cite[Proposition 3.8]{LimFineMWanti}. 
\epf

\section{Fine Selmer groups} \label{fine selmer sec}

Let $L$ be a finite extension of $F$ contained in $F_S$. Recall that the fine Selmer group of $E$ over $L$ is defined by
\[ R(E/L) =\ker\left(H^1(G_S(L),\Ep)\lra \bigoplus_{v\in S(L)} H^1(L_v, \Ep)\right).\]
Furthermore, the fine Selmer group and the classical Selmer group are related by the following exact sequence.

\bl
We have an exact sequence
\[ 0\lra R(E/L) \lra \Sel(E/L) \lra \bigoplus_{v\in S_p(L)}E(L_v)\ot_{\Zp}\Qp/\Zp.\]
In particular, the definition of the fine Selmer group does not depend on the
choice of the set $S$. \el

\bpf
 See \cite[Lemma 4.1]{LMu} or \cite[Section 2]{WuTS}.
\epf

Following Wuthrich \cite{WuTS}, the fine Mordell-Weil group $\mathcal{M}(E/L)$ is defined by
\[ \mathcal{M}(E/L) = \ker\Big(E(L)\ot_{\Zp}\Qp/\Zp \lra \bigoplus_{v\in S_p(L)} E(L_v)\ot_{\Zp}\Qp/\Zp \Big),\]
where $S_p(L)$ denotes the set of primes of $L$ above $p$. This fits into the following commutative diagram
\[   \entrymodifiers={!! <0pt, .8ex>+} \SelectTips{eu}{}\xymatrix{
    0 \ar[r]^{} & \mM(E/L) \ar[d] \ar[r] &  E(L)\ot_{\Zp}\Qp/\Zp
    \ar[d] \ar[r] & \displaystyle\bigoplus_{v\in S_p(L)} E(L_v)\ot_{\Zp}\Qp/\Zp  \ar@{=}[d]\\
    0 \ar[r]^{} & R(E/L) \ar[r]^{} & \Sel(E/L)\ar[r] & \displaystyle\bigoplus_{v\in S_p(L)} E(L_v)\ot_{\Zp}\Qp/\Zp
     } \]
with exact rows, where the leftmost vertical map is induced by the middle vertical map. Following Wuthrich \cite{WuTS}, the fine Tate-Shafarevich group $\Zhe(A/L)$ is given by
\[ \Zhe(E/L) = \coker\Big( \mM(E/L)\lra R(E/L)\Big).\]
Since the middle vertical map in the above diagram is injective, so is the leftmost vertical map. Consequently, a snake lemma argument yields a short exact sequence
\[ 0 \lra \mM(E/L) \lra R(E/L) \lra \Zhe(E/L)\lra 0\]
with $\Zhe(E/L)$ injecting into  $\sha(E/L)_{p^\infty}$, the $p$-primary part of the usual Tate-Shafarevich group.

\bigskip
We can now give the proof of our main theorem.

\bpf[Proof of Theorem \ref{mainthm}]
 Set $C_n$, resp., $C_\infty$, to be the cokernel of $\mathcal{M}(E/K_n)\lra E(K_n)\ot\Qp/\Zp$, resp., the cokernel of $\mathcal{M}(E/F_\infty)\lra E(F_\infty)\ot\Qp/\Zp$. On the other hand, from the argument in \cite[Proposition 3.14]{Van}, we see that 
\begin{equation}\label{MW H-rank}
  \rank_{\Zp\ps{H}} \big(E(F_\infty)\otimes\Qp/\Zp\big)^\vee= \left\{
  \begin{array}{ll}
    0, & \hbox{if $\epsilon(E/F,1)=+1$;} \\
    1, & \hbox{if $\epsilon(E/F,1)=-1$.}
  \end{array}
\right. 
\end{equation}
Plainly, if $\epsilon(E/F,1)=+1$, then $\mathcal{M}(E/F_\infty)$ has trivial $\Zp\ps{H}$-corank, and so by a result of Venjakob \cite{Ven} (or see \cite[Lemma 5.1]{LimFine}), $\mathcal{M}(E/F_\infty)^\vee$ is pseudo-null over $\Zp\ps{G}$. 

Now suppose that we have $\epsilon(E/F,1)=-1$. In view of (\ref{MW H-rank}) and the following tautological short exact sequence 
\[ 0 \lra \mathcal{M}(E/F_\infty)\lra E(F_\infty)\ot\Qp/\Zp \lra C_\infty \lra 0, \]
it therefore remains to show that $\rank_{\Zp\ps{H}} C_\infty >0$. 

Now consider the following commutative diagram
 \[   \entrymodifiers={!! <0pt, .8ex>+} \SelectTips{eu}{}\xymatrix{
    0 \ar[r] & \mathcal{M}(E/K_n)  \ar[d]^{\al_n} \ar[r] & E(K_n)\ot_{\Zp}\Qp/\Zp
    \ar[d]^{\be_n} \ar[r] & \bigoplus_{w\in S_p(K_n)} E(K_{n,w})\ot_{\Zp}\Qp/\Zp  \ar[d]^{\ga_n}\\
    0 \ar[r]^{} & \mathcal{M}(E/F_\infty)^{H_n} \ar[r]^{} & \Big(E(F_\infty)\ot_{\Zp}\Qp/\Zp\Big)^{H_n} \ar[r] &
    \left(\bigoplus_{u\in S_p(F_\infty)} E(F_{\infty, u})\ot_{\Zp}\Qp/\Zp \right)^{H_n}   } \]
    with exact rows.
Since $\ker \be_n$ is contained in $H^1(H_n, E(F_\infty)_{p^\infty})$, it follows from the argument in Proposition \ref{cyclotomic control} that the latter, and hence $\ker \be_n$, is finite. Similarly, we see that $\ker \ga_n$ is finite. Consequently, the map $C_n \lra (C_{\infty})^{H_n}$ also has a finite kernel.

We claim that $\mathrm{corank}_{\Zp} C_n \geq p^n-p^{n-1}$ for sufficiently large $n$.
Suppose that this claim holds. Since we have seen above that the map $C_n\lra (C_\infty)^H$ has finite kernel, it follows that $\mathrm{corank}_{\Zp}(C_\infty)^H\geq p^n-p^{n-1}$. On the other hand, since
\[ \mathrm{corank}_{\Zp}(C_\infty)^H = \Big(\mathrm{corank}_{\Zp\ps{H}} C_\infty \Big)p^n +O(1), \]
we therefore have 
\[ 
  \mathrm{corank}_{\Zp\ps{H}} C_\infty >0
\]
which is what we want to show.

It therefore remains to establish our claim. For this, we follow an idea of Lei \cite{LeiZ}. Consider the natural map
\[ f_n: E(K_n)\ot\Qp/\Zp \lra \bigoplus_{w\in S_p(K_n)} E(K_{n,w})\ot_{\Zp}\Qp/\Zp, \]
where $\mathrm{im}~ f_n =C_n$. Write $F^\mathrm{ac}$ for the anticyclotomic $\Zp$-extension of $F$, and $L_n$ the intermediate subextension of $F^\mathrm{ac}/F$ with $|L_n:F|=p^n$. Since $\epsilon(E/F,1)=-1$, we have
\[ \rank_\Z E(L_n) = p^n +O(1) \]
by \cite[Proposition 7.6]{Be}. As $L_n\subseteq K_n$, this in turn implies that  
$e_n \geq 1$ for sufficiently large $n$, where $e_n$ is defined as in Proposition \ref{MW control}. For each of such $n$, choose an element $x_n\in E(K_n)\setminus E(K_{n-1})$ which is of infinite order. Then for a prime $w$ of $K_n$ above $p$, the image of $x_n$ in $E(K_{n,w})$ is still of infinite order and so cannot be divisible by arbitrary power of $p$. Therefore, the image of $x_n$ in $E(K_{n,w})\ot\Qp/\Zp$ is nontrivial. 

As before, write $\Phi_n$ for the $p^n$th-cyclotomic polynomial viewed as an element in $\Zp\ps{H}$. For brevity, we write 
\[ V_n =T_p\Big(E(K_{n})\ot\Qp/\Zp\Big)\ot\Qp. \]
Then the map $f_n$ induces the following $\Qp[H/H_n]$-homomorphism
\[ V_n[\Phi_n] \hookrightarrow V_n \longrightarrow T_p\Big(E(K_{n,w})\ot\Qp/\Zp\Big)\ot\Qp \]
which is nonzero. Since $\Phi_n$ is irreducible over $\Qp[T]$, it follows that the image of the above homomorphism contains at least a copy of $\Qp[T]/\Phi_n$, and since $\mathrm{dim}_{\Qp}\big(\Qp[T]/\Phi_n\big) = p^n-p^{n-1}$, the said image has $\Qp$-dimension at least $p^n-p^{n-1}$. Consequently, one has $\mathrm{corank}_{\Zp}(C_n) \geq p^n-p^{n-1}$ for sufficiently large $n$. This establishes our claim and complete the proof of the theorem.
\epf

We end the paper with the proofs of the remaining two corollaries as stated in the introductory section.

\bpf[Proof of Corollary \ref{maincor1}]
For each $\Zp$-extension $L_\infty$ of $F$, we write $H_{L_\infty} = \Gal(F_\infty/L_\infty)$.
 Consider the following commutative diagram
 \[   \entrymodifiers={!! <0pt, .8ex>+} \SelectTips{eu}{}\xymatrix{
    0 \ar[r] &\mM(E/L_\infty)  \ar[d]^{\al} \ar[r] &  R(E/L_\infty)
    \ar[d]^{\delta}\\
    0 \ar[r]^{} & \mM(E/F_\infty)^{H_{L_\infty}} \ar[r]^{} & R(E/F_\infty)^{H_{L_\infty}} } \]
    with exact rows. From the proof of \cite[Proposition 3.8]{LimFineCam}, we see that the kernel of $\delta$ is cofinitely generated over $\Zp$ and hence so is the kernel of $\al$. On the other hand, combining Theorem \ref{mainthm} with \cite[Lemma 2.1(ii)]{LimFineCam}, we see that $\mM(E/F_\infty)^{H_{L_\infty}}$ is cotorsion over $\Zp\ps{\Gal(L_\infty/F)}$. Therefore, it follows that $\mM(E/L_\infty)$ is cotorsion over $\Zp\ps{\Gal(L_\infty/F)}$ as required.
\epf

\bpf[Proof of Corollary \ref{maincor2}]
Under the finiteness of $\Zhe(E/L)$, it is shown in \cite[Proposition 4.1]{LimFineDoc} that $\Zhe(E/L_\infty)$ is cotorsion over $\Zp\ps{\Gal(L_\infty/F)}$. The conclusion therefore follows from this, Corollary \ref{maincor1} and the following tautological short exact sequence
\[0\lra \mM(E/L_\infty) \lra R(E/L_\infty) \lra \Zhe(E/L_\infty) \lra 0.\]
\epf

\begin{acknowledgments}
Some part of the research of this article took place when the first named author
was visiting Harbin Engineering University and Harbin Institute of Technology, and
he would like to acknowledge the hospitality and conducive working conditions
provided by the said institutes. The second named author would especially like to
thank Chenyan Wu for her invaluable support and guidance during his visit to the
University of Melbourne. The second named author is also grateful to the University
of Melbourne for providing a stimulating research environment and the necessary
resources, which greatly contributed to the completion of this work. Finally, the second named author's research is supported by the National Natural Science Foundation of China under Grant No.
12001546, Heilongjiang Province under Grant No. 3236330122, and Harbin Engineering University
under Grant No. GK0000020127. The third named author's research is supported by the National Natural Science
Foundation of China under Grant No. 12331004.
\end{acknowledgments}

\end{document}